\documentclass[final,leqno,onefignum,onetabnum]{siamltex1213}
\usepackage{graphicx,epsfig,color}
\usepackage{booktabs}
\usepackage[notref,notcite]{showkeys}
\usepackage{amsmath,cases}
\usepackage{amssymb,bm,mathdots}
\usepackage{multirow}
\newtheorem{remark}[theorem]{Remark}
\input{siam.sty}

\title{Two spectral methods for 2D quasi-periodic scattering problems
} 

\author{Kui Du\thanks{
School of Mathematical Sciences and Fujian Provincial Key Laboratory of Mathematical Modelling and High-Performance Scientific Computation, Xiamen University, Xiamen 361005, China ({kuidu@xmu.edu.cn}). The research of this author was supported by the National Natural Science Foundation of China (No.11201392 and No.91430213), the Doctoral Fund of Ministry of Education of China (No.20120121120020), the Natural Science Foundation of Fujian Province of China (No.2013J01023) and the Fundamental Research Funds for the Central Universities (No.2013121003).} 
}

\begin{document}
\maketitle
\slugger{mms}{xxxx}{xx}{x}{x--x}

\begin{abstract} We consider the 2D quasi-periodic scattering problem in optics, which has been modelled by a boundary value problem governed by Helmholtz equation with transparent boundary conditions. A spectral collocation method and a tensor product spectral method are proposed to numerically solve the problem on rectangles. The discretization parameters can be adaptively chosen so that the numerical solution approximates the exact solution to a high accuracy. Our methods also apply to solve general partial differential equations in two space dimensions, one of which is periodic. Numerical examples are presented to illustrate the accuracy and efficiency of our methods. 
\end{abstract}



\pagestyle{myheadings}
\thispagestyle{plain}
\markboth{KUI DU}{TWO SPECTRAL METHODS FOR QUASI-PERIODIC SCATTERING}

\section{Introduction} With the development of electromagnetic and optics technology and the increasing demands in industrial and military applications, scattering from periodic structures  \cite{Joannopoulos2008photo} have attracted much interest in recent years. 
A variety of numerical methods including finite difference methods \cite{du2015finit}, finite element methods \cite{bao1995finit,bao2000numer,bao2000least,chen2003adapt, bao2005adapt,bao2006numer,wang2015adapt,du2015tenso}, spectral and spectral element methods  \cite{he2012effic,he2015spect}, integral equation methods  \cite{barnett2010new,barnett2011new,gillman2013fast,cho2015robust}, Dirichlet-to-Neumann map methods \cite{huang2006scatt}, and mode expansion method \cite{maes2005modeling,bai2007fourier} have been developed by the engineering community and the applied mathematical community for solving linear diffraction problems from periodic structures. 

Finite difference and finite element methods are easy to implement and result in sparse linear systems, which enable the use of sparse direct solvers. 
However, these methods typically have significant dispersion errors for high-frequency problems and thus fail to provide accurate solutions. If the medium is piecewise constant, boundary integral equations formulated on the interfaces of multilayer structures \cite{cho2015robust} or the boundaries of the obstacles \cite{barnett2011new,gillman2013fast} are natural and mathematically rigorous. By exploiting high-order quadratures, one obtains much higher efficiency and accuracy than finite difference and finite element methods. 
For problems with general medium, spectral methods \cite{fornberg1996pract,trefethen2000spect,boyd2011cheby,canuto2006spect,shen2011spect} can be employed. They are simple to implement and typically require relatively small number of unknowns to attain a fixed accuracy. We refer to \cite{townsend2014autom} for a brief survey of  existing spectral methods for partial differential equations.

In this paper, we consider 2D quasi-periodic scattering problems. By introducing transparent boundary conditions, the problem is governed by a Helmholtz equation with variable coefficients defined on rectangular domains. We propose a spectral collocation method and a tensor product spectral method. The first one uses the spectral collocation techniques \cite{fornberg1996pract,trefethen2000spect} and the second one, based on separable representations of the differential operator, combines the Fourier spectral method and the ultraspherical spectral method \cite{olver2013fast}. For vertically layered medium case, both methods only require to solve a one dimensional problem. The two spectral methods proposed in this paper can adaptively determine the number of unknowns and approximate the solution to a high accuracy. We remark that the tensor product spectral method leads to matrices that are banded or almost banded. By employing the fast algorithm in \cite{olver2013fast}, a layered medium approximation problem can be used as a preconditioner for the problem with general media.

Recently, Chebfun software \cite{driscoll2014chebf} was extended to solve problems in two dimensions \cite{townsend2013exten,townsend2014autom}, both of which are nonperiodic. We refer to \cite{trefethen2013appro} for the summary of the mathematics and algorithms of Chebyshev technology for nonperiodic functions. Extension to periodic functions in one dimension was proposed in \cite{wright2015exten}. Our spectral methods can be used to solve the problems in two space dimensions, one of which is periodic. Implementing our methods with the adaptive strategies of Chebfun software is easy.

The rest of this paper is organized as follows. In \S 2 the model problem is formulated and further reduced to a boundary value problem. In \S 3 we present the spectral collocation method  for the problem. Section 4 is devoted to the tensor product spectral method. Numerical examples illustrating the accuracy and efficiency of the methods are reported in \S 5. We present brief concluding remarks in \S 6.

\section{Problem setting} 
Assume that no currents are present and that the fields are source free. Then the electromagnetic fields in the whole space are governed by the following time-harmonic Maxwell equations (time dependence $\rme^{-\rmi\omega t}$) 
\beqa &&\nabla\times {\rm E}={\rmi\omega}{\mu}{\rm H}, \label{me1}\\ &&\nabla\times {\rm H}=-{\rmi\omega}\ve{\rm E},\label{me2}\eeqa
where $\rmi=\sqrt{-1}$ is the imaginary unit, $\omega$ is the angular frequency, $\mu$ is the permeability, $\ve$ is the permittivity, ${\rm E}$ is the electric field and $\rm H$ is the magnetic field.

In this paper, we consider a simple quasi-periodic scattering model: $\mu$ is constant everywhere; $\ve$ is periodic in $x$ variable of period $2\pi$, invariant in $z$ variable, and constant away from the region $|y|<1-\delta$, i.e., there exist constants $\ve^+$ and $\ve^-$ such that \beqas && \ve(x,y,z)=\ve^+, \mbox{ for } y> 1-\delta, \\ && \ve(x,y,z)=\ve^-, \mbox{ for } y< \delta-1. \eeqas  
In the TM polarization, the electric field ${\rm E}$ takes the simpler form
$${\rm E}(x, y, z) = {\rm E}(x, y) = (0, 0, u(x,y)).$$ The Maxwell equations (\ref{me1})-(\ref{me2}) yield the Helmholtz equation: \beq \Delta u+ \omega^2\ve\mu u=0.\label{reduv1}\eeq 

Consider the plane wave $u^{\rm i}=\rme^{\rmi\alpha_0 x-\rmi\beta_0 y}$ is incident from the above, where $\alpha_0=\omega\sqrt{\ve^+\mu}\sin\theta$, $\beta_0=\omega\sqrt{\ve^+\mu}\cos\theta$,  and $-\pi/2<\theta<\pi/2$ is the angle of incidence with respect to the positive $y$-axis. The incident wave $u^{\rm i}$ leads to reflected wave $u^{\rm r}$ and transmitted wave $u^{\rm t}$. For $y>1$, we have $u=u^{\rm i}+u^{\rm r}$, and for $y<-1$, $u=u^{\rm t}$. 
We are interested in quasi-periodic solution $u$ with phase $\rme^{\rmi\alpha_0 2\pi}$, i.e., $u\rme^{-\rmi\alpha_0 x}$ is $2\pi$-periodic in $x$ variable. Therefore, the reflected and transmitted waves can be written as
\beqas
&& u^{\rm r}=\sum_{j\in\mbbz} r_j\rme^{\rmi\alpha_jx+\rmi\beta_jy},\quad y>1,\\
&& u^{\rm t}=\sum_{j\in\mbbz} t_j\rme^{\rmi\alpha_jx-\rmi\gamma_jy}, \hspace{4mm} y<-1,
\eeqas where $r_j$ and $t_j$ are unknown complex scalar coefficients and \beqas&&\alpha_j=\alpha_0+j,\\ && \beta_j=\sqrt{\omega^2\ve^+\mu-\alpha_j^2},\qquad \im(\beta_j)\geq 0,\\ && \gamma_j=\sqrt{\omega^2\ve^-\mu-\alpha_j^2},\hspace{7.3mm} \im(\gamma_j)\geq 0.\eeqas
Note that if $\omega\sqrt{\ve^+\mu}$ is real, then the $j$'s satisfying $\beta_j>0$ correspond propagating modes. Throughout, we assume that $\beta_j\neq 0$ and $\gamma_j\neq 0$ for all $j\in\mbbz$. This assumption excludes the ``resonant'' cases, where waves can propagate along the $x$-axis. 

For a quasi-periodic function $f(x)$, define the linear operators $S$ and $T$ by \cite{bao1995mathe,bao1996numer,bao1997nonli}
\beqa && (Sf)(x)=\sum_{j\in\mbbz}\rmi\beta_j f_j\rme^{\rmi\alpha_j x},\label{mcals}\\ &&(T f )(x)=\sum_{j\in\mbbz}\rmi\gamma_j f_j \rme^{\rmi\alpha_j x},\label{mcalt}\eeqa where \beqas && f (x)=\sum_{j\in\mbbz}f_j \rme^{\rmi\alpha_j x},\\ && f_j =\frac{1}{2\pi}\int_0^{2\pi}f (x)\rme^{-\rmi\alpha_j x}\rmd x.\eeqas Let $\nu$ denote the unit outward normal. We obtain transparent boundary conditions: 
\beqa\label{dtn} && {\p_\nu u}-Su=-2\rmi\beta_0\rme^{\rmi\alpha_0 x-\rmi\beta_0},\qquad {\rm on}\quad (0,2\pi)\times \{1\}, \\&& {\p_\nu u}-T u=0,\hspace{27.1mm} {\rm on}\quad (0,2\pi)\times \{-1\}.\label{dtnn} \eeqa 
The quasi-periodic scattering problem is to solve the Helmholtz equation (\ref{reduv1}) in the rectangular domain $\O:=(0,2\pi)\times(-1,1)$ subject to the transparent boundary conditions (\ref{dtn})-(\ref{dtnn}) and the quasi-periodic boundary condition 
\beqa && u(x+2\pi,y)=\rme^{\rmi\alpha_0 2\pi}u(x,y).\label{qp}\eeqa

Define $v=u\rme^{-\rmi\alpha_0 x}$. Then $v$ satisfies \beq\label{periodic}\l\{\begin{array}{lll}\Delta_{\alpha_0} v+\omega^2\ve\mu v=0, &{\rm in} & \O=(0,2\pi)\times(-1,1),\\ {\p_y v}-\rme^{-\rmi\alpha_0 x}S (\rme^{\rmi\alpha_0 x}v)=-2\rmi\beta_0\rme^{-\rmi\beta_0}, &{\rm on} & (0,2\pi)\times{\{1\}}, \\ {\p_y v}+\rme^{-\rmi\alpha_0 x}T(\rme^{\rmi\alpha_0 x}v)=0, &{\rm on} & (0,2\pi)\times{\{-1\}},\\ v(x+2\pi,y)=v(x,y),& {\rm in} & \mbbr^2, \end{array}\r.\eeq where the operator $\Delta_{\alpha_0}$ is defined by $\Delta_{\alpha_0}=\p_{xx}+2\rmi\alpha_0\p_x-|\alpha_0|^2+\p_{yy}$. The solution to (\ref{periodic}) is unique at all but a discrete set of frequencies $\omega$ when the incident angle $\theta$ is fixed. Existence and uniqueness of the solution is strictly proved in Dobson \cite{dobson1993optim} by a variational approach. 
 
Since the operators $S$ and $T$ given in (\ref{mcals})-(\ref{mcalt}) are defined by infinite series, the computation has to be truncated in practice. For simplicity, we truncate $S$ and $T$ as follows:
\beqas && (S^N f)(x)=\sum_{j=1-q}^{N-q}\rmi\beta_jf_j\rme^{\rmi\alpha_j x}, \\ && (T^N f )(x)=\sum_{j=1-q}^{N-q}\rmi\gamma_j f_jT^N \rme^{\rmi\alpha_j x},\eeqas where $N$ is a sufficiently large integer and $q$ is an integer satisfying that $1\leq q\leq N$ and all propagating modes are contained in the middle of the truncated series. Next, we propose two spectral methods for the following problem  
\beq\label{tperiodic}\l\{\begin{array}{lll}\Delta_{\alpha_0} v+\omega^2\ve\mu v=0, &{\rm in} & \O=(0,2\pi)\times(-1,1),\\ {\p_y v}-\rme^{-\rmi\alpha_0 x}S^N(\rme^{\rmi\alpha_0 x}v)=-2\rmi\beta_0\rme^{-\rmi \beta_0}, &{\rm on} & (0,2\pi)\times{\{1\}}, \\ {\p_y v}+\rme^{-\rmi\alpha_0 x}T^N(\rme^{\rmi\alpha_0 x}v)=0, &{\rm on} & (0,2\pi)\times{\{-1\}},\\ v(x+2\pi,y)=v(x,y),& {\rm in} & \mbbr^2. \end{array}\r.\eeq
We refer to \cite[Section 3.1]{wang2015adapt} for the discussion on the existence and uniqueness of the solution of the problem (\ref{tperiodic}). In this paper, we focus on the accurate and efficient spectral methods for the problem (\ref{tperiodic}). Thus we always assume that the discrete problem has a unique solution. 

\section{A spectral collocation method}
We approximates the problem (\ref{tperiodic}) by Fourier discretization in $x$ variable and Chebyshev discretization in $y$ variable. 
Let $x_n=nh$ and $y_m=\cos(m\pi/M)$ with $h=2\pi/N$, $n=1,2,\ldots,N$, and $m=0,1,\ldots,M$. Introduce the first-order Fourier differentiation matrix \cite{gottlieb1984spect,weideman2000matla,trefethen2000spect}, ${\bf D}_x=[d_{ij}^x]_{i,j=1}^{N}$, \beqas  N \mbox{ odd },&\qquad& d_{ij}^x=\l\{\begin{array}{ll} 0, & i=j,\\ (-1)^{i-j}\dsp\frac{1}{2}\csc\frac{(i-j)h}{2}, & i\neq j, \end{array}\r.\\  N \mbox{ even }, &\qquad& d_{ij}^x=\l\{\begin{array}{ll} 0, & i=j,\\ (-1)^{i-j}\dsp\frac{1}{2}\cot\frac{(i-j)h}{2}, & i\neq j, \end{array}\r.\eeqas the second-order Fourier differentiation matrix \cite{gottlieb1984spect,weideman2000matla,trefethen2000spect}, ${\bf D}_{xx}=[d_{ij}^{xx}]_{i,j=1}^{N}$, \beqas N \mbox{ odd }, &\qquad& 
d_{ij}^{xx}=\l\{\begin{array}{ll} \dsp -\frac{\pi^2}{3h^2}+\frac{1}{12}, & i=j,\\ (-1)^{i-j+1}\dsp\frac{1}{2}\csc\frac{(i-j)h}{2}\cot\frac{(i-j)h}{2}, & i\neq j, \end{array}\r.\\  N \mbox{ even }, &\qquad& 
d_{ij}^{xx}=\l\{\begin{array}{ll} \dsp -\frac{\pi^2}{3h^2}-\frac{1}{6}, & i=j,\\ (-1)^{i-j+1}\dsp\frac{1}{2}\csc^2\frac{(i-j)h}{2}, & i\neq j, \end{array}\r.\eeqas and 
the $(M+1) \times (M+1)$ Chebyshev differentiation matrix \cite{canuto1988spect,weideman2000matla,trefethen2000spect}, ${\bf D}_y=[d_{ij}^y]_{i,j=0}^M,$  \beqas && d_{00}^y=\frac{2M^2+1}{6},\quad d_{MM}^y=-\frac{2M^2+1}{6},\\ && d_{ii}^y=-\frac{y_i}{2(1-y_i^2)},\quad i=1,\ldots,M-1,\\ && d_{ij}^y=\frac{c_i}{c_j}\frac{(-1)^{i+j}}{(y_i-y_j)}, \quad i\neq j, \quad i,j=0,1,\dots,M, \\ && c_i=\l\{\begin{array}{ll} 2, & i=0 \mbox{ or } M,\\ 1, & \mbox{otherwise.} \end{array}\r.\eeqas The discretization of the Helmholtz equation in (\ref{tperiodic}) takes the form \beq\label{discrete} {\bf P}{\bf V}{\bf D}_{xx}^\rmt+2\rmi\alpha_0{\bf P}{\bf V}{\bf D}_x^\rmt-|\alpha_0|^2{\bf P}{\bf V}+{\bf P}{\bf D}_{y}^2{\bf V}+\omega^2\mu {\bf P}({\bm \ve}\odot{\bf V})={\bf 0},\eeq where ${\bf P}$ is the $(M-1)\times(M+1)$ matrix obtained by deleting the first and last rows of the identity matrix ${\bf I}_{M+1}$, ${\bf V}=[v_{mn}]_{m=0,n=1}^{M,N}$, $v_{mn}$ is the approximate solution at the point $(x_n,y_m)$, ${\bm \ve}=[\ve_{mn}]_{m=0,n=1}^{M,N}$, $\ve_{mn}=\ve(x_n,y_m)$, and $\odot$ denotes the componentwise product.  
The operators $\rme^{-\rmi\alpha_0 x}S^N(\rme^{\rmi\alpha_0 x}\bm\cdot)$ and $\rme^{-\rmi\alpha_0 x}T^N(\rme^{\rmi\alpha_0 x}\bm\cdot)$ can be approximated by $N\times N$ matrices:
\beqa\label{S} && {\bf S}={\bf G}^*{\bm\Lambda}_\beta{\bf G}:={\bf G}^*\diag\{\rmi\beta_{1-q},\rmi\beta_{2-q},\cdots,\rmi\beta_{N-q}\}{\bf G},\\ \label{T}&& {\bf T}={\bf G}^*{\bm\Lambda}_\gamma{\bf G}:={\bf G}^*\diag\{\rmi\gamma_{1-q},\rmi\gamma_{2-q},\cdots,\rmi\gamma_{N-q}\}{\bf G},\eeqa where ${\bf G}$ is the $N\times N$ matrix with $(i,j)$ entry $\rme^{-\rmi 2(i-q)j\pi/N }/\sqrt{N}$, $i,j=1,\ldots,N$. The discrete forms of the transparent boundary conditions are given by \beqa \label{beta1} && {\bf e}_0^\rmt{\bf D}_y{\bf V}-{\bf e}_0^\rmt{\bf V}{\bf S}^\rmt=-2\rmi\beta_0\rme^{-\rmi\beta_0} {\bf e}^\rmt, \\ \label{gamma1}&& {\bf e}_M^\rmt{\bf D}_y{\bf V}+{\bf e}_M^\rmt{\bf V}{\bf T}^\rmt = {\bf 0},\eeqa where \beqas&& {\bf e}_0=[1 \ 0 \ \cdots \ 0]^\rmt\in\mbbr^{M+1},\\ && {\bf e}_M=[0 \ \cdots \ 0\ 1]^\rmt\in\mbbr^{M+1},\\&& {\bf e}=[1\ 1\ \cdots \ 1]^\rmt\in\mbbr^{N}.\eeqas 

Let 
\beq\label{helm}{\bf H}= ({\bf D}_{xx}+2\rmi\alpha_0{\bf D}_x-|\alpha_0|^2 {\bf I}_N)\otimes {\bf P}+{\bf I}_N\otimes({\bf PD}_y^2)+\omega^2\mu({\bf I}_N\otimes {\bf P}) \diag\{\vec({\bm\ve})\}.\eeq Combining (\ref{discrete}), (\ref{beta1}) and (\ref{gamma1}) yields the linear system \beq\label{global} {\bf A}{\bf v}={\bf g},\eeq where \beqs {\bf A}=\l[\begin{array}{c}{\bf I}_N\otimes({\bf e}_0^\rmt{\bf D}_y)-{\bf S}\otimes {\bf e}_0^\rmt\\ {\bf H} \\ {\bf I}_N\otimes({\bf e}_M^\rmt{\bf D}_y)+{\bf T}\otimes {\bf e}_M^\rmt\end{array}\r],\eeqs and 
$${\bf v}=\vec({\bf V}),\qquad
{\bf g}=\l[\begin{array}{c}-2\rmi\beta_0\rme^{-\rmi\beta_0} {\bf e}\\ {\bf 0}\\ {\bf 0}\end{array}\r]\in\mbbc^{MN+N}.$$ 

\begin{remark} The matrix $\bf P$ in {\rm(\ref{discrete})} is called the downsampling matrix {\rm\cite{driscoll2015recta}} obtained by using the inner Chebyshev points of the second kind,$$y_m=\cos\l(\frac{m\pi}{M}\r),\qquad m=1,2,\ldots,M-1,$$ as resampling points. One might choose other points, such as the Chebyshev points of the first kind, $$y_m=\cos\l(\frac{(2m-1)\pi}{2(M-1)}\r),\qquad m=1,2,\ldots, M-1.$$ We refer to {\rm\cite{xu2015expli}} for explicit construction of rectangular differentiation matrix ${\bf PD}_y^2$.  
\end{remark}
  
\subsection{Layered medium case} Assume that the medium in $\O$ is vertically layered, i.e., $\ve(x,y)=\ve(y)$. Then, the matrix ${\bf H}$ in (\ref{helm}) takes the form $${\bf H}= ({\bf D}_{xx}+2\rmi\alpha_0{\bf D}_x-|\alpha_0|^2 {\bf I}_N)\otimes {\bf P}+{\bf I}_N\otimes({\bf P}({\bf D}_y^2+\omega^2\mu{\bm\ve}^y)),$$ where $${\bm\ve}^y=\diag\{\ve(y_0),\ve(y_1),\cdots,\ve(y_M)\}.$$ 
Let ${\bf F}$ be the   
discrete Fourier transform matrix with $(i,j)$ entry $\rme^{-\rmi2(i-1)(j-1)\pi/N}/\sqrt{N}$,  $i,j=1,2,\ldots,N.$ We have the following propositions. The proofs are straightforward.

\begin{proposition}\label{pfft} The matrices ${\bf S}$ in {\rm(\ref{S})} and ${\bf T}$ in {\rm(\ref{T})} satisfy $${\bf S}={\bf F}^*{\bm\Lambda}_{S}{\bf F},\qquad{\bf T}={\bf F}^*{\bm\Lambda}_{T}{\bf F},$$ where \beqas{\bm\Lambda}_{S}&=&\diag\{\lambda_1^\beta,\lambda_2^\beta,\cdots,\lambda_{N}^\beta\}\\ &=&\diag\{\rmi\beta_0,\rmi\beta_1,\cdots,\rmi\beta_{N-q},\rmi\beta_{1-q},\rmi\beta_{2-q},\cdots,\rmi\beta_{-1}\},\eeqas and \beqas{\bm\Lambda}_{T}&=&\diag\{\lambda_1^\gamma,\lambda_2^\gamma,\cdots,\lambda_{N}^\gamma\}\\ &=&\diag\{\rmi\gamma_0,\rmi\gamma_1,\cdots,\rmi\gamma_{N-q},\rmi\gamma_{1-q},\rmi\gamma_{2-q},\cdots,\rmi\gamma_{-1}\}.\eeqas  
\end{proposition}  \vspace{-5mm}
\begin{proposition}\label{diag} The first-order spectral differentiation matrix ${\bf D}_x$ and the second-order spectral differentiation matrix ${\bf D}_{xx}$ can be diagonalized by ${\bf F}$. Specifically, we have $${\bf D}_x={\bf F}^*{\bm\Lambda}_x{\bf F}, \qquad {\bf D}_{xx}={\bf F}^*{\bm\Lambda}_{xx}{\bf F},$$ where \beqas{\bm \Lambda}_x&=&\diag\{\lambda_1^x,\lambda_2^x,\cdots,\lambda_{N}^x\}\\&=&\l\{\begin{array}{ll}\rmi\cdot\diag\{0,1,\cdots,(N-1)/2,(1-N)/2,(3-N)/2,\cdots,-1\}, & \mbox{if N is odd}, \\ \rmi\cdot\diag\{0,1,\cdots,N/2-1,0,1-N/2,2-N/2,\cdots,-1\}, & \mbox{if N is even}, \end{array}\r.\eeqas and \beqas{\bm \Lambda}_{xx}&=&\diag\{\lambda_1^{xx},\lambda_2^{xx},\cdots,\lambda_{N}^{xx}\}\\&=&\l\{\begin{array}{ll}-\diag\{0,1,4,\cdots,((N-1)/2)^2,((1-N)/2)^2,\cdots,4,1\}, & \mbox{if N is odd}, \\ -\diag\{0,1,4,\cdots,(N/2)^2,(N/2-1)^2,\cdots,4,1\}, & \mbox{if N is even}. \end{array}\r.\eeqas
\end{proposition}

Multiplying the linear system (\ref{global}) by $$\l[\begin{array}{ccc} {\bf F} &&\\ &{\bf F}\otimes {\bf I}_{M-1} &\\ &&{\bf F}\end{array}\r]$$ from the left yields 
\beq\label{uncoupled}\l[\begin{array}{c}{\bf I}_N\otimes({\bf e}_0^\rmt{\bf D}_y)-{\bm\Lambda}_{S}\otimes {\bf e}_0^\rmt\\ ({\bm\Lambda}_{xx}+2\rmi\alpha_0{\bm\Lambda}_x-|\alpha_0|^2 {\bf I}_N)\otimes {\bf P}+{\bf I}_N\otimes({\bf P}({\bf D}_y^2+\omega^2\mu{\bm\ve}^y)) \\ {\bf I}_N\otimes({\bf e}_M^\rmt{\bf D}_y)+{\bm\Lambda}_{T}\otimes {\bf e}_M^\rmt\end{array}\r]\wh{\bf v}=\wh{\bf g},\eeq where $$\wh{\bf v}=({\bf F}\otimes {\bf I}_{M+1}){\bf v},\qquad \wh{\bf g}=\l[\begin{array}{c}-2\rmi\beta_0\rme^{-\rmi\beta_0} \wh{\bf e}\\ {\bf 0}\\ {\bf 0}\end{array}\r],\qquad \wh{\bf e}={\bf Fe}=\l[\begin{array}{c}\sqrt{N}\\ {\bf 0}\end{array}\r].$$  
Reordering the unknowns and equations in (\ref{uncoupled}), we obtain $N$ uncoupled linear systems of order $M+1$: \beq\label{oned}\l[\begin{array}{c}{\bf e}_0^\rmt {\bf D}_y-\rmi\beta_0{\bf e}_0^\rmt \\{\bf P}({\bf D}_y^2+\omega^2\mu{\bm\ve}^y)-|\alpha_0|^2{\bf P}\\ {\bf e}_M^\rmt {\bf D}_y+\rmi\gamma_0{\bf e}_M^\rmt\end{array}\r]\wh{\bf v}_1=\l[\begin{array}{c} -2\rmi\beta_0\rme^{-\rmi\beta_0}\sqrt{N}\\ {\bf 0}\\ 0\end{array}\r],\eeq and for $i=2,\ldots,N$, $$\l[\begin{array}{c}{\bf e}_0^\rmt {\bf D}_y-\lambda_i^\beta{\bf e}_0^\rmt \\(\lambda_i^{xx}+2\rmi\alpha_0\lambda_i^x-|\alpha_0|^2){\bf P}+{\bf P}({\bf D}_y^2+\omega^2\mu{\bm\ve}^y)\\ {\bf e}_M^\rmt {\bf D}_y+\lambda_i^\gamma{\bf e}_M^\rmt\end{array}\r]\wh{\bf v}_i={\bf 0},$$ where $\wh{\bf v}_i$ is the $i$-th column of ${\bf VF}^\rmt$. Obviously, $\wh{\bf v}_i={\bf 0}$  for $i\neq 1$. After solving the linear system (\ref{oned}), we obtain ${\bf v}$ by ${\bf v}=({\bf F}^*\otimes{\bf I}_{M+1})\wh{\bf v}$.

\section{A tensor product spectral method} We propose a tensor product spectral method (see \S 4.3) for the problem (\ref{tperiodic}) by combining the Fourier spectral method and the ultraspherical spectral method \cite{olver2013fast}.

\subsection{Fourier spectral method for periodic second order linear ordinary differential equations}
We consider the second-order linear ordinary differential equation (ODE) \beq\label{pode} \mcall w:=a(x)w''(x)+b(x)w'(x)+c(x)w(x)=f(x),\qquad x\in[0,2\pi] \eeq with periodic boundary conditions. Here, $a(x)$, $b(x)$, $c(x)$ and $f(x)$ are periodic functions on $[0,2\pi]$. The Fourier spectral method finds an infinite vector $${\bf w}=\l[\begin{array}{ccccccc} \cdots & w_{-2} & w_{-1} & w_{0} & w_{1} & w_{2} & \cdots  \end{array}\r]^\rmt,$$ such that the Fourier expansion of the solution of (\ref{pode}) is given by $$w(x)=\sum_{j\in\mbbz}w_j \rme^{\rmi j x},\qquad x\in[0,2\pi].$$

Note that we have $$w'(x)=\sum_{j\in\mbbz}\rmi j w_j \rme^{\rmi jx},\qquad w''(x)=\sum_{j\in\mbbz} (-j^2)w_j \rme^{\rmi jx}.$$ The first-order differentiation operator is given by $$\mcald_x=\rmi\l[\begin{array}{ccccccccc} \ddots&&&&&&&& \\ &-3&&&&&&& \\ &&-2&&&&&& \\&&&-1&&&&& \\ &&&&0&&&&\\ &&&&&1&&&\\&&&&&&2&&\\ &&&&&&&3& \\&&&&&&&&\ddots \end{array}\r],$$ and the second-order differentiation operator is given by $\mcald_x^2.$ 
In order to handle variable coefficients of the forms $a(x)w''(x)$, $b(x)w'(x)$ and $c(x)w(x)$ in (\ref{pode}), we need to represent the multiplication of two Fourier series as an operator on coefficients. Let $\mcalt[a]$ denote the multiplication operator that represents multiplication of two Fourier series, i.e., if $\bf w$ is a vector of Fourier expansion coefficients of $w(x)$, then $\mcalt[a]{\bf w}$ returns the Fourier expansion coefficients of $a(x)w(x)$. Suppose that $a(x)$ is given by its Fourier series $$a(x)=\sum_{j\in\mbbz} a_j \rme^{\rmi jx}.$$ Then the explicit formula for $\mcalt[a]$ is given by the following Toeplitz operator $$\mcalt[a]=\l[\begin{array}{ccc} \ddots & \vdots & \ddots\\ \ddots & a_{-2} & \ddots\\ \ddots & a_{-1} & \ddots\\ \ddots& a_0 & \ddots \\ \ddots & a_1 & \ddots \\ \ddots & a_2 & \ddots \\ \ddots& \vdots & \ddots \end{array}\r].$$
This multiplication operator looks dense; however, if $a(x)$ is approximated by a {\it trigonometric polynomial of degree $n$}, then $\mcalt[a]$ is banded with a bandwidth of $n$. 
Combining the differentiation and multiplication operators yields $$(\mcalt[a]\mcald_x^2+\mcalt[b]\mcald_x+\mcalt[c]){\bf w}={\bf f},$$ where ${\bf w}$ and ${\bf f}$ are vectors of Fourier expansion coefficients of $w(x)$ and $f(x)$, respectively.
We need truncate the operator to derive a practical numerical scheme. 
Let $\mcalp_N$ be the projection operator satisfying $$\mcalp_N{\bf w}=\l[\begin{array}{ccc} {\bf 0} & {\bf I}_N & {\bf 0}  \end{array}\r]{\bf w}=\l[\begin{array}{cccc} w_{1-q} & w_{2-q} & \cdots & w_{N-q} \end{array}\r]^\rmt.$$ We obtain the following linear system 
$$(\mcalp_N\mcall^\mbbd\mcalp_N^\rmt)(\mcalp_N{\bf w})=\mcalp_N{\bf f},$$ where $$\mcall^\mbbd=\mcalt[a]\mcald_x^2+\mcalt[b]\mcald_x+\mcalt[c].$$  The truncation parameter $N$ can be adaptively chosen so that the numerical solution approximates the exact solution to relative machine precision.

\subsection{The ultraspherical spectral method \cite{olver2013fast} for second order linear ordinary differential equations}
Consider the second order linear ordinary differential equation \beq\label{code}\mcall w:=a(x)w''(x)+b(x)w'(x)+c(x)w(x)=f(x),\eeq where $a(x)$, $b(x)$, $c(x)$, $w(x)$ and $f(x)$ are functions defined on $[-1,1]$. The ultraspherical spectral method finds an infinite vector $${\bf w}=\l[\begin{array}{ccc}w_0 & w_1 & \cdots\end{array}\r]^\rmt$$ such that the Chebyshev expansion of the solution of (\ref{code}) is given by $$w(x)=\sum_{j=0}^\infty w_j T_j(x),\qquad x\in[-1,1],$$ where $T_j(x)$ is the degree $j$ Chebyshev polynomial. 

Note that we have the following recurrence relation $$\frac{\rmd^\lambda T_n}{\rmd x^\lambda}=\l\{\begin{array}{ll} 2^{\lambda-1}n(\lambda-1)!C_{n-\lambda}^{(\lambda)}, & n\geq\lambda,\\ 0, & 0\leq n\leq \lambda-1, \end{array}\r.$$ where $C_j^{(\lambda)}$ is the ultraspherical polynomial with an integer parameter $\lambda\geq 1$ of degree $j$ \cite{olver2010nist}.
Then the differentiation operator for the $\lambda$th derivative is given by \beq\label{l0}\mcald_\lambda=2^{\lambda-1}(\lambda-1)!\l[\begin{array}{ccccc} {\bf 0}&\lambda&&& \\ &&\lambda+1&&\\ &&&\lambda+2&\\ &&&&\ddots\end{array}\r], \quad \lambda\geq 1.\eeq Here $\bf 0$ in (\ref{l0})
denotes a $\lambda$-dimensional zero row vector. The conversion operator converting a vector of Chebyshev expansion coefficients to a vector of $C^{(1)}$ expansion coefficients, denoted by $\mcals_0$, and the conversion operator converting a vector of $C^{(\lambda)}$ expansion coefficients to a vector of $C^{(\lambda+1)}$ expansion coefficients, denoted by $\mcals_\lambda$, are given by $$\mcals_0=\l[\begin{array}{ccccc}1 &0 &-\frac{1}{2}&& \\ &\frac{1}{2}&0&-\frac{1}{2}&\\ &&\frac{1}{2}&0 &\ddots\\ &&&\frac{1}{2}&\ddots\\ &&&&\ddots \end{array}\r],\quad \mcals_\lambda=\l[\begin{array}{ccccc}1 &0 &-\frac{\lambda}{\lambda+2}&& \\ &\frac{\lambda}{\lambda+1}&0&-\frac{\lambda}{\lambda+3}&\\ &&\frac{\lambda}{\lambda+2}&0 &\ddots\\ &&&\frac{\lambda}{\lambda+3}&\ddots\\ &&&&\ddots \end{array}\r],\ \lambda\geq 1.$$ We also require the multiplication operator $\mcalm_0[a]$ that represents multiplication of two Chebyshev series, and the multiplication operator $\mcalm_\lambda[a]$ that represents multiplication of two $C^{(\lambda)}$ series, i.e., if ${\bf w}$ is a vector of Chebyshev expansion coefficients of $w(x)$, then $\mcalm_0[a]{\bf w}$ returns the Chebyshev expansion coefficients of $a(x)w(x)$, and  $\mcalm_\lambda[a]\mcals_{\lambda-1}\cdots\mcals_0{\bf w}$ returns the $C^{(\lambda)}$ expansion coefficients of $a(x)w(x)$. Suppose that $a(x)$ has the Chebyshev expansion $$a(x)=\sum_{j=0}^\infty a_jT_j(x).$$ Then $\mcalm_0[a]$ can be written as  \cite{olver2013fast}:
$$\mcalm_0[a]=\frac{1}{2}\l[\begin{array}{ccccc}2a_0 & a_1 & a_2 & a_3 & \cdots\\ a_1 &2a_0 & a_1 & a_2 &\ddots \\ a_2 & a_1 & 2a_0 & a_1 & \ddots \\ a_3 & a_2& a_1& 2a_0& \ddots \\ \vdots & \ddots &\ddots &\ddots &\ddots \end{array}\r]+\frac{1}{2}\l[\begin{array}{ccccc} 0 & 0 & 0 & 0 & \cdots\\ a_1 &a_2 & a_3 & a_4 &\cdots \\ a_2 & a_3 & a_4 & a_5 & \iddots \\ a_3 & a_4& a_5& a_6& \iddots \\ \vdots & \iddots &\iddots &\iddots &\iddots \end{array}\r].$$
The explicit formula for the entries of $\mcalm_\lambda[a]$ with $\lambda\geq 1$ is given in \cite{olver2013fast}. This multiplication operators $\mcalm_\lambda[a]$ with $\lambda\geq 0$  look dense; however, if $a(x)$ is approximated by a truncation of its Chebyshev or $C^{(\lambda)}$ series, then $\mcalm_\lambda[a]$ is banded. Combining the differentiation, conversion and multiplication operators yields \beq\label{dcode}(\mcalm_2[a]\mcald_2+\mcals_1\mcalm_1[b]\mcald_1+\mcals_1\mcals_0\mcalm_0[c]){\bf w}=\mcals_1\mcals_0{\bf f},\eeq where ${\bf w}$ and ${\bf f}$ are vectors of Chebyshev expansion coefficients of $w(x)$ and $f(x)$, respectively.  We need truncate the operator to derive a practical numerical scheme. Let $\mcalq_M$ be the projection operator given by $$\mcalq_M=\l[\begin{array}{cc}  {\bf I}_M & {\bf 0}  \end{array}\r].$$ We obtain the following linear system $$(\mcalq_M\mcall^\mbbd\mcalq_M^\rmt)(\mcalq_M{\bf w})=(\mcalq_M\mcals_1\mcals_0\mcalq_M^\rmt)(\mcalq_M{\bf f}),$$ where $$\mcall^\mbbd=\mcalm_2[a]\mcald_2+\mcals_1\mcalm_1[b]\mcald_1+\mcals_1\mcals_0\mcalm_0[c].$$

\subsection{Discretization of (\ref{tperiodic}) for the special case  $\ve(x,y)=\phi(x)\psi(y)$} We seek the approximate solution $\tilde v(x,y)$ to the problem (\ref{tperiodic}), $$\tilde v(x,y)= \sum_{i=0}^{M-1}\sum_{j=1-q}^{N-q}v_{ij}\rme^{\rmi j x}T_i(y).$$ The linear differential operator in (\ref{tperiodic}) can be expressed as $$\mcall^x\otimes \mcali+\mcali\otimes\mcall^y+\omega^2\mu\phi\otimes\psi,$$ where $\mcall^x:=\p_{xx}+2\rmi\alpha_0\p_x-|\alpha_0|^2$, $\mcall^y:=\p_{yy}$, $\mcali$ is the identity operator, and $\phi$ and $\psi$ are the corresponding multiplication operations. The discretization of the Helmholtz equation in (\ref{tperiodic}) takes the form
\beq\label{tpfc} {\bf Q}{\bf C}{\bf V}{\bf X}^\rmt+{\bf Q}{\bf Y}{\bf V}+\omega^2\mu {\bf Q}{\bm \Psi}{\bf V}{\bm \Phi}^\rmt={\bf 0},\eeq where {\bf Q} is the $(M-2)\times M$ matrix obtained by deleting the last two rows of the identity matrix ${\bf I}_M$, ${\bf C}=\mcalq_M\mcals_1\mcals_0\mcalq_M^\rmt$, ${\bf X}=\mcalp_N(\mcald_x^2+2\rmi\alpha_0\mcald_x)\mcalp_N^\rmt-|\alpha_0|^2{\bf I}_N$, ${\bf Y}=\mcalq_M\mcald_2\mcalq_M^\rmt$, ${\bm \Phi}=\mcalp_N\mcalt[\phi]\mcalp_N^\rmt$, ${\bm \Psi}=\mcalq_M\mcals_1\mcals_0\mcalm_0[\psi]\mcalq_M^\rmt$,  ${\bf V}=[v_{ij}]_{i=0,j=1-q}^{M-1,N-q}$.  The discrete forms of the transparent boundary conditions are given by 
\beqa \label{beta2} && {\bf b}_1^\rmt{\bf V}-{\bf a}_1^\rmt{\bf V}{\bm\Lambda}_\beta=-2\rmi\beta_0\rme^{-\rmi\beta_0} {\bf e}_q^\rmt, \\ \label{gamma2}&& {\bf b}_2^\rmt{\bf V}+{\bf a}_2^\rmt{\bf V}{\bm\Lambda}_\gamma = {\bf 0},\eeqa where \beqas&& {\bm\Lambda}_\beta=\diag\{\rmi\beta_{1-q},\rmi\beta_{2-q},\cdots,\rmi\beta_{N-q}\},\\&& {\bm\Lambda}_\gamma=\diag\{\rmi\gamma_{1-q},\rmi\gamma_{2-q},\cdots,\rmi\gamma_{N-q}\},\\&& {\bf b}_1=[0 \ 1 \ 4 \ \cdots \ (M-1)^2]^\rmt,\\ && {\bf b}_{2}=[0 \ 1 \ -4 \ \cdots \ (-1)^M(M-1)^2]^\rmt, \\&& {\bf a}_1=[1\ 1\ \cdots \ 1]^\rmt\in\mbbr^{M},\\&& {\bf a}_2=[1\ -1\ \cdots \ (-1)^{(M-1)}]^\rmt\in\mbbr^{M},\eeqas and ${\bf e}_q$ is the $q$-th column of the identity matrix ${\bf I}_N$. 
 
Let
\beq\label{helmfc}{\bf H}= {\bf X}\otimes ({\bf QC})+{\bf I}_N\otimes({\bf QY})+\omega^2\mu{\bm \Phi}\otimes ({\bf Q}{\bm \Psi}).\eeq 
 Combining (\ref{tpfc}), (\ref{beta2}) and (\ref{gamma2}) yields the linear system \beq\label{globalfc} {\bf A}{\bf v}={\bf g},\eeq where \beqs {\bf A}=\l[\begin{array}{c}{\bf I}_N\otimes{\bf b}_1^\rmt-{\bm\Lambda}_\beta\otimes {\bf a}_1^\rmt\\ {\bf I}_N\otimes{\bf b}_2^\rmt+{\bm\Lambda}_\gamma\otimes {\bf a}_2^\rmt\\ {\bf H} \end{array}\r],\eeqs and $${\bf v}=\vec({\bf V}),\qquad
{\bf g}=\l[\begin{array}{c}-2\rmi\beta_0\rme^{-\rmi\beta_0} {\bf e}_q\\ {\bf 0}\\ {\bf 0}\end{array}\r]\in\mbbc^{MN}.$$ 

\begin{remark} For general bivariate function $\ve(x,y)$, we can use its low rank approximant, i.e., sum of functions of the form $\phi(x)\psi(y)$, where $\phi(x)$ and $\psi(y)$ are univariate functions. An algorithm {\rm \cite{townsend2013exten,townsend2013gauss}} which is mathematically equivalent to Gaussian elimination with complete pivoting can be used to construct low rank approximations.  
\end{remark}

\subsubsection{Layered medium case} In this case, we have $\ve(x,y)=\psi(y)$, i.e., $\phi(x)=1$. Then $\bm\Phi={\bf I}_N.$  The matrix ${\bf H}$ in (\ref{helmfc}) takes the form $${\bf H}= {\bf X}\otimes ({\bf QC})+{\bf I}_N\otimes({\bf QY}+\omega^2\mu{\bf Q}{\bm \Psi}).$$  
Reordering the unknowns and equations in (\ref{globalfc}), we obtain  \beq\label{onedfc}\l[\begin{array}{c}{\bf b}_1^\rmt -\rmi\beta_0{\bf a}_1^\rmt \\{\bf b}_2^\rmt +\rmi\gamma_0{\bf a}_2^\rmt\\ {\bf Q}({\bf Y}+\omega^2\mu{\bm\Psi}-|\alpha_0|^2{\bf C})\end{array}\r]{\bf v}_q=\l[\begin{array}{c} -2\rmi\beta_0\rme^{-\rmi\beta_0}\\ 0\\ {\bf 0}\end{array}\r],\eeq and \beq\label{ofc}\l[\begin{array}{c}{\bf b}_1^\rmt -\rmi\beta_{j-q}{\bf a}_1^\rmt \\{\bf b}_2^\rmt +\rmi\gamma_{j-q}{\bf a}_2^\rmt\\ {\bf Q}({\bf Y}+\omega^2\mu{\bm\Psi}-((j-q)^2+2\alpha_0(j-q)+|\alpha_0|^2){\bf C})\end{array}\r]{\bf v}_j={\bf 0},\qquad j\neq q,\eeq  where ${\bf v}_j$ is the $j$-th column of ${\bf V}$. Obviously, ${\bf v}_j={\bf 0}$  for $j\neq q$.

Note that the fast algorithm in \cite{olver2013fast} can be used to solve the linear system (\ref{onedfc}), which requires $\mcalo(m^2M)$ operations. Here, $m$ is the number of Chebyshev points needed to resolve the function $\psi(y)$. The coefficient matrix resulting from  
a layered medium approximation problem can be used as a preconditioner for the problem with general medium. The corresponding computational complexity for the preconditioner solve is $\mcalo(m^2MN)$.

\section{Numerical results} We have performed numerical experiments in numerous cases. In this section, we present a few typical results of these experiments to illustrate the accuracy and efficiency of the two spectral methods. All computations are performed with MATLAB R2012a.

\begin{figure}[!htpb]
\centerline{\epsfig{figure=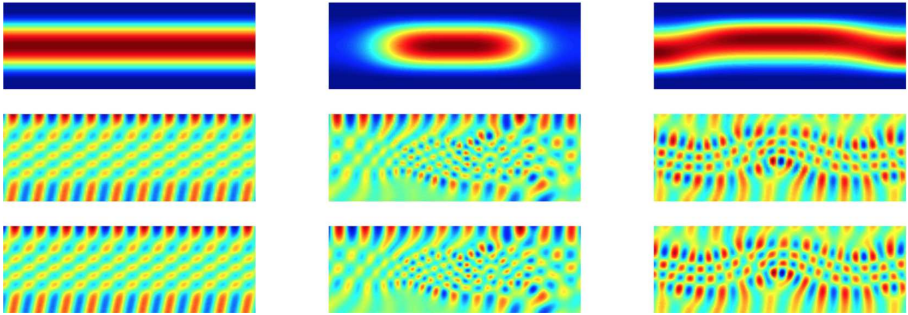,height=2in}}
\caption{Media {\rm(Top)} and corresponding real parts of the waves obtained by the spectral collocation method {\rm(Middle)} and by the tensor product spectral method {\rm(Bottom)}. From left to right{\rm:} numerical results for $\ve_1(x,y)$, $\ve_2(x,y)$, and $\ve_3(x,y)$, respectively.}\label{nr}
\end{figure}

The parameters are chosen as $\theta=3\pi/7$, $\omega=10$, $\ve^+=\ve^-=1$, $\mu=1$. We consider three cases: \beqas
&&\ve_1(x,y)=\ve(y)=1+\dsp\rme^{\frac{3}{y^2-1}+4},\\
&&\ve_2(x,y)=1+\dsp\rme^{\frac{3}{y^2-1}+4-\cos(\pi\sin(x/2))},\\
&&\ve_3(x,y)=1+\dsp\rme^{\frac{3}{y^2-1}+4-y\cos(\pi\sin(x/2))}. 
\eeqas
The medium characterized by $\ve_1(x,y)$ is vertically layered. The bivariate function $\ve_2(x,y)$ is of rank $2$. We use the rank $7$ approximant obtained by the algorithm proposed in {\rm \cite{townsend2013exten,townsend2013gauss}} to approximate $\ve_3(x,y)$. In Figure \ref{nr}, we plot the three media and the real parts of the corresponding waves obtained by our spectral methods. We observe that the numerical results obtained by the two methods coincide well for all the three media. 

\begin{figure}[!htpb]
\centerline{\epsfig{figure=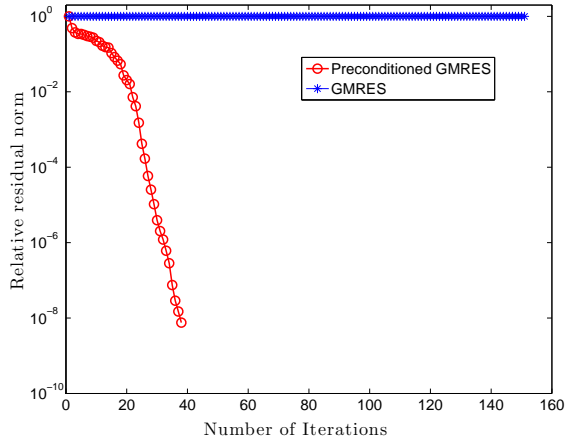,height=2.5in}}
\caption{Convergence history for {\rm GMRES} applied to the linear system  {\rm(\ref{globalfc})}.}\label{conhis}
\end{figure}

We also test the performance of the problem with $\ve_1(x,y)$ as a preconditioner for the problem with $\ve_3(x,y)$. The fast algorithm in \cite{olver2013fast} is used to solve the linear systems (\ref{onedfc}) and (\ref{ofc}). The GMRES algorithm \cite{saad1986gmres} is used as the iterative solver. The initial guess is set to be the zero vector. GMRES with the layered medium preconditioner obtained a solution with the relative residual norm less than $10^{-8}$ at the $37$th iteration, while GMRES without preconditioning almost stagnates; see Figure \ref{conhis} for the convergence history.

%

\section{Concluding remarks} We have proposed a spectral collocation method and a tensor product spectral method for solving the 2D quasi-periodic scattering problem. Both of the methods can adaptively determine the number of unknowns and approximate the solution to a high accuracy. The tensor product spectral method is more interesting because it leads to matrices that are banded or almost banded, which enable the use of the fast, stable direct solver \cite{olver2013fast}. Based on this fast solver, a layered medium preconditioning technique is used to solve the problem with general media. Our methods also apply to solve general partial differential equations in two space dimensions, one of which is periodic. Extension of our methods to bi-periodic structure diffraction grating problem \cite{bao1997varia} is being considered.


\end{document}